\documentclass[a4paper,11pt]{article}

\usepackage[latin1]{inputenc}
\usepackage{fancyhdr}
\usepackage{eurosym,amssymb,amsmath,verbatim,amsthm}
\usepackage{lastpage,full page,color}
\usepackage{setspace}
\usepackage{hyperref}
\usepackage{tabularx}
\usepackage[section]{placeins}

\usepackage{setspace}

\usepackage{multicol}
\usepackage{multirow}
\usepackage{color}
\usepackage{colortbl}
\usepackage{xcolor}

\theoremstyle{plain}

\newtheorem*{theorem*}{Theorem}
\newtheorem*{mtheorem*}{Main Theorem}

\usepackage{hyperref}
\hypersetup{
    colorlinks=true,
    linkcolor=blue,
    filecolor=magenta,      
    urlcolor=cyan,
}

  \renewenvironment{thebibliography}[1]{%
    \begin{oldthebibliography}{#1}%
      \setlength{\parskip}{0ex}%
      \setlength{\itemsep}{0ex}%
  }%
  {%
    \end{oldthebibliography}%
  }

\begin{document}
\author{J.~Schillewaert}
\title{Solution to Bishnoi's conjecture on minimal $t$-fold blocking sets of maximal size}
%Hjelmslev planes and Segre Varieties as variants of Hermitian Veroneseans
\date{}
\maketitle

\begin{abstract}
Bishnoi conjectured that if a minimal $t$-fold blocking set in a projective plane of prime power order has maximal size then it is either a projective plane minus one point, the complement of a Baer subplane or a unital. In this note we prove this conjecture.
\end{abstract}

%\section{Introduction}

A $t$-fold blocking set in a projective plane $\pi$ of order $n$ is a set of points $\mathcal{K}$ such that each line of $\pi$ intersects $\mathcal{K}$ in at least $t$ points and some line of $\pi$ intersects $\mathcal{K}$ in exactly $t$ points. It is called minimal if every point of $\mathcal{K}$ is contained in a line intersecting $\mathcal{K}$ in exactly $t$ points.

The case $t=1$ is well-studied going back to Bruen and Thas \cite{BT}. They have shown that in a finite projective plane of order $n$ the size of a (1-fold) minimal blocking set is bounded above by $n\sqrt{n}+1$. This bound is sharp in the case when $n$ is a square, and $S$ is a unital in $\pi$, i.e. a set of $n\sqrt{n}+1$ points in $\pi$ such that each line of $\pi$ intersects $S$ in 1 or $\sqrt{n}+1$ points. An example of a unital in a Desarguesian plane of square order $\mathrm{PG}(2,q^2)$ is the point set $\mathcal{K}$ of a non-degenerate Hermitian variety, i.e. the points $(x,y,z) \in \mathrm{PG}(2,q^2)$ satisfying $x^{q+1}+y^{q+1}+z^{q+1}=0$.

Multiple blocking sets were introduced by Bruen in \cite{Bruen} and lower bounds were obtained by Ball in \cite{Ball}.
In \cite{AB} Bishnoi proved the first general upper bound on the size of minimal $t$-fold blocking sets extending the result of Bruen and Thas to the case of general $t$. More precisely he proved the following

\begin{theorem*} [\cite{AB}]\label{Bishnoimain}
A minimal $t$-fold blocking set $S$ in a finite projective plane of order $n$ has size at most

$$\frac{1}{2}n\sqrt{4tn-(3t+1)(t-1)}+\frac{1}{2}(t-1)n+t$$
If the size of $S$ is equal to this upper bound, then every line of $\pi$ intersects $S$ in exactly $t$ or 
$b+1 := \frac{1}{2}(\sqrt{4tn-(3t+1)(t-1)}+t-1)+1$ points.
\end{theorem*}

In particular for the case $t=1$ he recovers the result of Bruen and Thas, for $t=n$ the bound is equal to $n^2+n$, and $S$ is a projective plane minus a point, and when $n$ is a square for $t=n-\sqrt{n}$ we obtain the upper bound $n^2-\sqrt{n}$, in which case $S$ is the complement of a Baer subplane. A Baer subplane $\mathcal{B}$ in $\pi$ is a set of $n+\sqrt{n}+1$ points of $\pi$ such that each line of $\pi$ intersects $\mathcal{B}$ in 1 or $\sqrt{n}+1$ points.  
For example, in a Desarguesian plane of square order $\mathrm{PG}(2,q^2)$ the subplane $\mathrm{PG}(2,q)$ is a Baer subplane.

\begin{mtheorem*}[Bishnoi's conjecture] In case $n$ is a prime power $q$ the only possible values for $t$ for which equality can be reached in Theorem \ref{Bishnoimain} are 
\begin{itemize}
\item $t=1$ when $q$ is a square, in this case $S$ is a unital in $\pi$. 
\item $t=q-\sqrt{q}$ when $q$ is a square, in this case $S$ is the complement of a Baer subplane in $\pi$.
\item $t=q$ for any $q$, in this case $S$ is the plane $\pi$ with one point removed.
\end{itemize} 
\end{mtheorem*}

\begin{proof}
From the expression for $b$ we obtain $b^2+b(1-t)-t+t^2=tn$ $(\dagger)$.
Also $b-t+1$ divides $n$ $(\star)$, see e.g. Section 2 of \cite{PR}.
Assume now that $n=q=p^k$, then by $(\star)$ we have $b-t+1=p^h$. Write $t=\alpha p^l$, with $(\alpha,p)=1$.
 Hence $(\dagger)$ becomes
\begin{equation} \label{eqn}
p^h(p^h+\alpha p^l-1)-\alpha p^l+\alpha^2 p^{2l} = \alpha p^{l}q
\end{equation}

By $(\dagger)$ $t$ is a divisor of $b(b+1)=(p^h-1+t)(p^h+t)$, hence $t=\alpha p^l$ is a divisor of $p^h(p^h-1)$, so $l\leq h$ and $\alpha$ divides $p^h-1$ $(\ast)$. We will distinguish four cases.

{\bf Case I} $l>0,h>0,l<h$: dividing (\ref{eqn}) by $p^l$ implies that $p$ divides $\alpha$, a contradiction.

{\bf Case II} $l>0,h>0,l=h$: Then after division of (\ref{eqn}) by $p^h$ we obtain
\begin{equation}\label{eqn2}
(p^h+\alpha p^h-1)-\alpha+\alpha^2 p^{h} = \alpha q
\end{equation}
So $\alpha = \beta p^h-1$ and by (*) $\alpha = p^h-1$, and (\ref{eqn2}) yields after simplification that $p^{2h}=q$, hence {\color{blue} $q$ is a square, $t=q-\sqrt{q}$ and $b=q-1$, and $S$ is the complement of a Baer subplane}.

{\bf Case III} $l=0,h>0$: Then from (\ref{eqn}) we get $p^h(p^h+\alpha-1)-\alpha+\alpha^2 = \alpha q$ which implies that $\alpha^2-\alpha$ is a multiple of $p^h$. Since $(\alpha,p)=1$ we must have $\alpha = \beta p^h+1$. By $(\ast)$ $\beta=0$, so $\alpha=1$ implying $p^{2h}=q$ and hence {\color{blue}  $q$ is a square, $t=1$ and $b=\sqrt{q}$, and $S$ is a unital}.

{\bf Case IV} $l=h=0$: Then we have {\color{blue} $b=t=q$, and $S$ is the projective plane with one point removed}.

\end{proof}

\end{document}